\documentclass[11pt]{article}

\usepackage{amsfonts,amsmath,amssymb}
\usepackage{mathrsfs,mathtools,url}
\usepackage{latexsym}
\usepackage{tikz,color}
\usepackage[T1]{fontenc}
\usepackage{array}
\usepackage{xspace}
\usepackage{cite}
\usepackage{multirow}
\usepackage{graphicx}
\usepackage{pgfplots}

\usepackage[T1]{fontenc}
\usepackage{xcolor} 
\usepackage{txfonts} 
\usepackage{hyperref}

\usetikzlibrary{arrows}

\newtheorem{theorem}{Theorem}[section]

\newtheorem{lemma}[theorem]{Lemma}
\newtheorem{proposition}[theorem]{Proposition}
\newtheorem{observation}[theorem]{Observation}

\newtheorem{corollary}[theorem]{Corollary}
\newtheorem{remark}[theorem]{Remark}

\newtheorem{problem}[theorem]{Problem}

\DeclareMathOperator {\gp} {gp}
\DeclareMathOperator {\diam} {diam}

\let\deg\relax
\DeclareMathOperator {\deg} {deg}
\DeclareMathOperator {\mut} {\mu_{\rm t}}
\DeclareMathOperator {\muo} {\mu_{\rm o}}
\DeclareMathOperator {\mud} {\mu_{\rm d}}

\newcommand{\proof}{\noindent{\bf Proof.\ }}
\newcommand{\qed}{\hfill $\square$ \bigskip}
\newcommand{\smallqed}{{\tiny ($\Box$)}}

\title{Visibility in graphs under edge and vertex removal}

\author{Pakanun Dokyeesun$^{a}$, \enskip  Csilla Bujt\'{a}s$^{a, b}$\\\\
	$^{a}$ \small Institute of Mathematics, Physics and Mechanics, Ljubljana, Slovenia \\
	\small {\tt papakanun@gmail.com \quad ORCID: 0009-0004-4705-2542}\\
	$^{b}$ \small Faculty of Mathematics and Physics, University of Ljubljana, Slovenia\\
	\small {\tt csilla.bujtas@fmf.uni-lj.si\quad ORCID: 0000-0002-0511-5291}\\
}

\date{\empty}

\begin{document}
	
	\maketitle
	
	\begin{abstract}
		For a connected graph $G$ and $X\subseteq V(G)$, we say that two vertices $u$, $v$ are $X$-visible if there is a shortest $u,v$-path $P$ with $V(P)\cap X \subseteq \{u,v\}$. If every two vertices from $X$ are $X$-visible, then $X$ is a mutual-visibility set in $G$. The largest cardinality of such a set in $G$ is the mutual-visibility number $\mu(G)$. When the visibility constraint is extended to further types of vertex pairs, we get the definitions of outer, dual, and total mutual-visibility sets and the respective graph invariants $\muo(G)$, $\mud(G)$, and $\mut(G)$.
		
		This work concentrates on the possible changes in the four visibility invariants when an edge $e$ or a vertex $x$ is removed from $G$ and the graph remains connected. It is proved that $\frac{1}{2}\mu(G) \le \mu(G-e) \le 2\mu(G)$\,  and\, 
		$\frac{1}{6}\muo(G) \le \muo(G-e) \le 2\muo(G)+1$ hold for every graph. Further general upper bounds established here are $\mut(G-e) \leq \mut(G)+2$\, and\, $\mu(G-x) \leq 2\mu(G)$. For all but one of the remaining cases, it is shown that the visibility invariant may increase or decrease arbitrarily under the considered local operation. For example, neither $\mud(G-e)$ nor $\mud(G-x)$ allows lower or upper bounds of the form $a \cdot \mud(G)+b$ with a positive constant $a$. Along the way, the realizability of the four visibility invariants in terms of the order is also characterized in the paper. 
	\end{abstract}
	
	
	\medskip\noindent
	\textbf{AMS Math.\ Subj.\ Class.\ (2020)}: 05C12, 05C31, 05C69
	
	
	\section{Introduction} \label{sec:intro}
	
	For a connected graph $G$ and a vertex set $X \subseteq V(G)$, two vertices \emph{$u$ and $v$ are $X$-visible} if there is a shortest $u,v$-path $P$ such that $V(P) \cap X  \subseteq \{u,v\}$.
	We define four properties and graph invariants based on the visibility of vertices. Let $X$  $\subseteq V(G)$.
	\begin{description}
		\item[$X$ is a mutual-visibility set] in $G$ if every two vertices $u,v$ with $u,v \in X$ are $X$-visible. The maximum cardinality of such a set is the \emph{mutual-visibility number} $\mu(G)$.
		\item[$X$ is an outer mutual-visibility set] in $G$ if every two vertices $u,v$ with $u \in X$ and $v \in V(G)$ are $X$-visible. The maximum cardinality of such a set is the \emph{outer mutual-visibility number} $\muo(G)$.
		\item[$X$ is a dual mutual-visibility set] in $G$ if every two vertices $u,v$ with $u,v \in X$ or with $u,v \in V(G)\setminus X$ are $X$-visible. The maximum cardinality of such a set is the \emph{dual mutual-visibility number} $\mud(G)$.
		\item[$X$ is a total mutual-visibility set] in $G$ if every two vertices $u,v \in V(G)$ are $X$-visible. The maximum cardinality of such a set is the \emph{total mutual-visibility number} $\mut(G)$.      
	\end{description}
	
	\paragraph{Earlier studies.} Di Stefano~\cite{distefano-2022} initiated the study of mutual-visibility sets in graphs motivated by well-studied problems related to distributed systems and robot navigation. This work was soon followed by further studies and the introduction of total mutual-visibility~\cite{cicerone-2023+}, outer and dual mutual-visibility sets~\cite{CiDiDrHeKlYe-2023}. These visibility sets and invariants have been extensively investigated in recent years, see for example~\cite{haxenovich-2024a+, boruzanli-2024, BresarYero-2024, Bujtas-2025+, cicerone-2025, cicerone-2023a, cicerone-2024b, roy-2025, tian-2024}.
	
	From among the different directions of the research, we highlight the works on visibility in hypercubes and Hamming graphs that show that even good estimations for the visibility invariants are hard to prove on these graph classes~\cite{haxenovich-2024a+, Bujtas, korze-2024, korze-2025, kuziak-2023}. It was also pointed out that Zarankiewicz's celebrated problem, the covering problem from design theory, specific Tur\'an-type problems on graphs and hypergraphs have close interrelation with visibility problems on some graph classes (see e.g.,~\cite{boruzanli-2024, BresarYero-2024, Bujtas, cicerone-2023, cicerone-2024b}).
	
	\paragraph{Our results.}
	This manuscript studies the changes in the visibility invariants of a graph when a vertex $x$ or an edge $e$ is removed from $G$ and the graph remains connected. As it turns out, the question is more challenging than expected at first sight. For every  $\sigma \in \{\mu, \muo, \mud, \mut\}$, we concentrate on general lower and upper bounds on $\sigma(G-x)$ and $\sigma(G-e)$ of the form $a\cdot \sigma(G) +b$, where $a,b$ are fixed constants and $a>0$. A general bound of the given form or its non-existence is proved for all but one case. 
	Our results are summarized in Table~\ref{table:result}. Proofs for edge removal and vertex removal are presented in Sections~\ref{sec:edge-removal} and~\ref{sec:vertex-removal}, respectively.

	\setlength{\tabcolsep}{16pt}
	\renewcommand{\arraystretch}{2}
	\begin{table}[h]
		\centering
		\begin{tabular}{ |c||c |c |c |  }
			\hline
			Name& Lower bound & Upper bound & Reference\\
			\hline
			\hline
			$\mu(G-e)$  & $\frac{1}{2}\, \mu(G)$   & $2\,\mu(G)$ & Theorem \ref{thm:edge-mu}\\
			\hline
			$\muo(G-e)$ & $\frac{1}{6}\, \muo(G)$  & $2\muo(G)+1$  & Theorem \ref{thm:edge-outer} \\
			\hline
			$\mud(G-e)$ & NO & NO & Proposition \ref{prop:edge-mud} \\
			\hline
			$\mut(G-e)$ & NO   & $\mut(G) + 2$ & Theorem \ref{thm:edge-mut} \\
			\hline
			$\mu(G-x)$  & NO & $2\,\mu(G)$ &  Corollary~\ref{cor:lower-vertex}, Theorem~\ref{thm:mu-vertex}\\
			\hline
			$\muo(G-x)$ & NO & ?  &  Corollary~\ref{cor:lower-vertex}\\
			\hline
			$\mud(G-x)$ & NO & NO &  Corollary~\ref{cor:lower-vertex}, Corollary~\ref{cor:upper-vertex}\\
			\hline
			$\mut(G-x)$ & NO & NO  &  Corollary~\ref{cor:lower-vertex}, Corollary~\ref{cor:upper-vertex}\\
			\hline
		\end{tabular}
		\caption{Summary of our results. ``NO'' means that no general bound of the required form $a\cdot \sigma(G) +b$ exists. The question mark indicates the case that remains open.} 
		\label{table:result}
	\end{table}
	Along the way of studying visibility invariants under local operators, the following problem arose and has been solved. Let $p$ and $q$ be nonnegative integers and $\sigma \in \{\mu, \muo, \mud, \mut\}$. We say that the pair $(p,q)$ is $(\sigma, n)$-realizable, if there exists a graph $G$ with $\sigma(G)=p$ and $n(G)=q$. In Section~\ref{sec:realizability}, we characterize $(\sigma, n)$-realizable pairs for each $\sigma \in \{\mu, \muo, \mud, \mut\}$.
	
	
	\section{Preliminaries and terminology}
	\label{sec:pre}
	\subsection{Standard definitions}
	
	In this paper, all graphs $G=(V(G),E(G))$ are simple and connected. The \emph{order} of the graph is $n(G)=|V(G)|$.
	The {\em degree} $\deg(u)$ of a vertex $u \in V(G)$ is the number of vertices adjacent to $u$ in $G$. The {\em maximum degree} $\Delta(G)$ of $G$ is $\max \{\deg(u): u \in V(G)\}$. 
	A {\em leaf} is a vertex of degree one. 
	We say that $u$ is a \emph{universal vertex} if $\deg(u)=n(G)-1$.  The {\em open neighborhood } $N_G(u)$ of a vertex $u$ is the set of neighbors of $u$ in $G$ and the {\em closed neighborhood }  $N_G[u] = N_G(u)\cup \{u\}$. 
	If $A\subseteq V(G)$, the subgraph of $G$ induced by $A$ is denoted by $G[A]$. As usual, a path, cycle, and complete graph of order $n$ are denoted by $P_n$, $C_n$, and $K_n$, respectively. An \emph{internal vertex} of a path $P$ is a vertex of degree $2$ in $P$. For a positive integer $k$, notation $[k]$ stands for the set $\{1,\ldots, k\}$.
	
	The {\em distance} $d_G(u,v)$ between vertices $u$ and $v$ of $G$ is the number of edges on a shortest $u,v$-path. The {\em diameter} $\diam(G)$ of $G$ is the maximum distance between pairs of vertices of $G$. Let $H$ be a subgraph of $G$. Then $H$ is {\em isometric} if for each pair of vertices $u,v\in V(H)$ we have $d_H(u,v) = d_G(u,v)$, and $H$ is {\em convex} if for any vertices $u,v\in V(H)$, every shortest $u,v$-path in $G$ lies in $H$. 
	
	A set $X \subseteq V(G)$ is a \emph{$\mu$-set} in $G$, if $X$ is a mutual-visibility set of maximum cardinality in $G$. We use the terms \emph{$\muo$-set}, \emph{$\mud$-set}, and \emph{$\mut$-set} with analogous meanings. When an edge $e \in E(G)$ is removed from the graph $G$, we delete $e$ from the edge set and obtain the graph $G-e$. When a vertex $x \in V(G)$ is removed, we delete $x$ and all incident edges to obtain $G-x$. Note that $G-x$ is exactly $G[V(G)\setminus \{x\}]$.
	
	\subsection{Lemmas from the literature}
	
	By definition, the following inequalities hold for every connected graph $G$:
	\begin{equation} \label{eq:1}
		\mut(G) \leq  \muo(G) \leq \mu(G) \quad \mbox{and} \quad \mut(G) \leq  \mud(G) \leq \mu(G)  
	\end{equation}
	
	It is easy to show that the subsets of a mutual-visibility set are also mutual-visibility sets in $G$. The analogous property holds for outer and total mutual-visibility sets.
	\begin{lemma}
		{\rm \cite{CiDiDrHeKlYe-2023}}
		\label{lem:subset-closed}
		If $X$ is a mutual-visibility set (resp.\ outer, total mutual-visibility set) of a graph $G$ and $Y\subseteq X$, then $Y$ is a mutual-visibility set (resp.\ outer, total mutual-visibility set) of $G$.
	\end{lemma}
	For dual mutual-visibility sets, the situation is different. For example, in a $5$-cycle, two neighbors $v_1, v_2$ form a $\mud$-set, but no one-element subset is a dual mutual-visibility set. In~\cite{Bujtas-2025+}, among more general results, a graph $F_t$ is constructed for every $t \ge 2$ such that it admits exactly one $t$-element $\mud$-set but no non-empty subset of it is a dual mutual-visibility set in $F_t$.
	The following type of monotonicity, however, holds for all varieties of visibility sets.
	\begin{lemma}{\rm \cite{Bujtas-2025+, distefano-2022}}
		\label{lem:convex}
		Let $X$ be a mutual-visibility set (dual, outer, total mutual-visibility set) of a graph $G$. If $H$ is a convex subgraph of $G$, then $X\cap V(H)$ is a mutual-visibility set (dual, outer, total mutual-visibility set) of $H$.
	\end{lemma}
	By definition, a total mutual-visibility set of a path $P_3$ never contains the middle vertex. By Lemma~\ref{lem:convex}, if $v$ is the middle vertex of a convex $P_3$-subgraph of $G$, then $v$ does not belong to any total mutual-visibility set of $G$. 
	As the following statement shows, we may prove that a set $X \subseteq V(G)$ is a total mutual-visibility set in $G$ by checking the visibility of vertices $u,v \in V(G)$ with $d_G(u,v)=2$. 
	\begin{lemma}{\rm \cite{Bujtas-2025+}}
		\label{lem:distance two}
		If $G$ is a connected graph and $X\subseteq V(G)$, then $X$ is a total mutual-visibility set of $G$ if and only if every two vertices $u, v \in V(G)$ with $d_G(u,v) = 2$ are $X$-visible.  
	\end{lemma} 
	
	For a tree $T$, all four visibility invariants have the same value, which is equal to the number of leaves in $T$. For a cycle $C_n$, the following equalities hold.
	\begin{proposition}{\rm \cite{CiDiDrHeKlYe-2023,  distefano-2022, tian-2024}} 
		\label{prop:cycles}
		For every $n \ge 3$, cycle $C_n$ satisfies the following statements.
		\begin{itemize}
			\item[(i)] $\mu(C_n)=3$;
			\item[(ii)] $\muo(C_3)=3$ and $\muo(C_n)=2$ if $n \ge 4$;
			\item[(iii)] $\mud(C_3)=\mud(C_4)=3$,  $\mud(C_5)=\mud(C_6)=2$, and $\mud(C_n)=0$ if $n \ge 7$;
			\item[(iv)] $\mut(C_3)=3$, $\mut(C_4)=2$, and  $\mut(C_n)=0$ if $n \ge 5$.
		\end{itemize}    
	\end{proposition}
	
	Let $G$ be a graph and $X\subseteq V(G)$. $X$ is 
	a {\em general position set} of $G$ if $V(P)\cap X = \{u,v\}$ for every pair $u,v \in X$ and every shortest $u,v$-path $P$. The {\em general position number} $\gp(G)$ of $G$ is the largest cardinality of such a set in $G$. It follows directly from definitions that $\mu(G) \ge \gp(G)$. For more details, see~\cite{chandran-2016,  Klavzar-2025+, manuel-2018, Thomas-2024a}.
	The possible changes in the general position number of a graph under vertex and edge removal were studied in~\cite{dokyeesun-2025}. Since invariants $\mu(G)$ and $\gp(G)$ are closely related, our proof methods are sometimes similar to those in~\cite{dokyeesun-2025}. 
	
	\subsection{Specific preliminaries}
	
	In this subsection, we present two lemmas and introduce some terminology with basic observations that will be used later in several proofs. 
	
	\begin{lemma} \label{lem:mud-1}
		If $G$ is a graph with $\mud(G)=1$ and $\{x\}$ is a $\mud$-set in $G$, then $\mut(G-x)=0$.
	\end{lemma}
	\proof Let $G$ be a graph with a $\mud$-set $X=\{x\}$. Since every two vertices from $V(G)\setminus  X$ are $X$-visible, $G-x$ is an isometric subgraph in $G$. Suppose now that $\mut(G-x) \ge 1$. Then, by Lemma~\ref{lem:subset-closed}, there exists a vertex $y$ in $G-x$ such that $\{y\}$ is a total mutual-visibility set in $G-x$. We prove that $Y=\{x,y\}$ is a dual mutual-visibility set in $G$. Vertices $x$ and $y$ are clearly $Y$-visible. If $u$ and $v$ are two vertices from $V(G)\setminus Y$ then, since $\{y\}$ is a total mutual-visibility set in $G-x$,  there is a shortest $u,v$-path $P$ in $G-x$ which does not go through $y$. Since $G-x$ is an isometric subgraph, $P$ remains a shortest $u,v$-path in $G$. Hence, $u$ and $v$ are $Y$-visible in $G$. Then $Y$ is a dual mutual-visibility set in $G$, which contradicts the condition $\mud(G)=1$. Therefore, $\mut(G-x)=0$ as stated. \qed
	
	\begin{lemma}
		\label{lem:lower}
		Let $G$ be a connected graph. If there is a partition of $V(G)$ into $A$ and $B$ such that $G[A]$ is a complete graph $K_k$ and $G[B]$ is an isometric subgraph of $G$, then  $\mud(G) \ge k$. 
	\end{lemma}
	\proof
	Since $G[A]$ is complete, $A$ is a mutual-visibility set of $G$ by definition. 
	Suppose now that $A$ is not a dual mutual-visibility set of $G$. Then there exist $u,v \in B$ which are not $A$-visible. Thus, every shortest $u,v$-path in $G$ must contain a vertex from $A$ that implies $d_G(u,v) < d_{G[B]}(u,v)$. It contradicts the condition of $G[B]$ being isometric. Therefore, $A$ is a dual mutual-visibility set of $G$, and we conclude  $\mud(G) \ge k$.
	\qed
	
	When we apply Lemma~\ref{lem:lower}, $G[B]$ is frequently an induced subgraph of diameter $2$, which is always isometric. Further, if $G$ contains a cycle $C_n$, for $n \in \{3,4,5,6\}$, and two adjacent vertices $v_1,v_2 \in V(C_n)$ with $\deg_G(v_1)=\deg_G(v_2)=2$, then $G-\{v_1,v_2\}$ is an isometric subgraph and hence, $\mud(G) \ge 2$.

	\paragraph{Introducing $Z_G(x,y)$.} For a graph $G$, a set $Z \subseteq V(G)$ and two vertices $x,y \in V(G)$, we define 
	$$ Z_G(x,y)= \{u: d_G(u,x)\leq d_G(u,y), \enskip u \in Z\}.
	$$
	If $Z$ is a mutual-visibility set (resp.\ an outer or a total mutual-visibility set) in $G$, then Lemma~\ref{lem:subset-closed} implies that $Z_G(x,y)$ and $Z \setminus Z_G(x,y)$ are also mutual-visibility sets (resp.\ outer or total mutual-visibility sets) in $G$. Let $u,u' \in Z_G(x,y)$. By definition of $Z_G(x,y)$ and since a shortest path between two vertices is always an isometric subgraph, we deduce that no shortest $u,u'$-path contains both $x$ and $y$. (Otherwise, one vertex from $\{u,u'\} \subseteq Z_G(x,y)$ would be closer to $y$ than to $x$, a contradiction.)
	Further, under the same condition $u,u' \in Z_G(x,y)$, if there is a shortest $u,u'$-path through $y$, then $d_G(u,x)= d_G(u,y)$ and $d_G(u',x) = d_G(u',y)$ hold and therefore, a shortest $u,u'$-path through $x$ also exists.  By definition, $Z_G(x,y) \cup Z_G(y,x) =Z$ holds and consequently, $\max\{|Z_G(x,y)|, |Z_G(y,x)|\} \ge |Z|/2$. 
	
	If $e=xy$ is an edge in the graph $G$, then for every $Z \subseteq V(G)$ and $u \in Z_G(x,y)$, no shortest $u,x$-path goes through $y$. In $G-e$, all these paths remain shortest paths and $d_{G-e}(u,x)= d_G(u,x)$. Since the removal of $e=xy$ cannot decrease any distances in the graph, we infer that $d_{G-e}(u,x) \leq d_{G-e}(u,y)$. Thus $u \in Z_G(x,y)$ implies $u \in Z_{G-e}(x,y)$. The proof is similar for the other direction, and we obtain $Z_G(x,y)= Z_{G-e}(x,y)$ if $e=xy$.
	
	\paragraph{Removing blocking vertices.} Let $Z$ be a mutual-visibility set in a graph $G$ and $w \in V(G)$. For every $u \in Z$, select a shortest $w,u$-path $P_{w,u}$ that contains the least internal vertices from $Z$; that is $|V(P_{w,u})\cap Z|$ is minimum. Since $Z$ is a mutual-visibility set, a shortest $w,u$-path with two internal vertices $u_1$ and $u_2$ from $Z$ can be replaced by a shortest $w,u$-path with only one internal vertex from $Z$. Indeed, if the path is $u\dots u_1 \dots u_2 \dots w$, we can replace the part $u \dots u_1 \dots u_2$ with a shortest $u,u_2$-path with no internal vertex from $Z$.
	
	When $w \in V(G)$ and a shortest path $P_{w,u}$ is fixed for every $u \in Z$, we say that an ordered pair $(u',u)$ with $u',u \in Z$ and $u' \neq u$ is a \emph{blocking pair} or that \emph{$u'$ blocks $u$}, if $P_{w,u}$ contains $u'$. The set $B_{w,Z}$ of blocking vertices contains a vertex $u' \in Z$ if $u'$ blocks at least one $u\in Z$. As $P_{w,u}$ is fixed, every $u \in Z$ is blocked by at most one vertex. Further, if $(u',u)$ is a blocking pair, then no blocking pair $(z,u')$ or $(u,z)$ exists. As follows, we have $|B_{Z,w}| \leq |Z|/2 $ for every $w \in V(G)$. We also note that after removing the blocking vertices from $Z$, every vertex $u$ in $Z\setminus B_{Z,w}$ will see $w$ via the path $P_{w,u}$. When we refer to a \emph{blocking set} $B_{Z,w}$ it is always assumed that a set of shortest paths $P_{w,u}$ is fixed.
	
	\medskip
	
	The main properties explained above are summarized in the following observation. 
	\begin{observation} \label{obs:Z-and-blocking}
		Let $Z$ be a mutual-visibility set in a graph $G$ and $x,y,w \in V(G)$. Let $B_{Z,w}$ be a blocking set in $G$.
		\begin{itemize}
			\item[(i)] If $u,u' \in Z_G(x,y)$, then no shortest $u,u'$-path contains both $x$ and $y$. Further, if a shortest $u,u'$-path contains $y$, then there is a shortest $u,u'$-path containing $x$.
			\item[(ii)] $\max\{|Z_G(x,y)|, |Z_G(y,x)|\} \ge |Z|/2$.
			\item[(iii)] If $e=xy$ is an edge in $G$, then $Z_G(x,y)= Z_{G-e}(x,y)$. 
			\item[(iv)] If $Z_w= Z\setminus B_{Z,w}$, then  $|Z_{w}| \ge |Z|/2$ and, for every vertex $u \in Z_{w}$, vertices $u$ and $w$ are $Z_w$-visible.
		\end{itemize}    
	\end{observation}

	\section{Edge removal}  \label{sec:edge-removal}
	This section aims to analyze the potential variations in the four visibility invariants resulting from the removal of edges. In the first subsection, we present observations on a family of graphs. Subsequently, we consider each visibility invariant and determine how the removal of an edge may affect its value.
	
	\subsection{A useful example}
	We open this subsection by defining a family of graphs. The changes in their visibility invariants under edge removal will give (sharp) examples for our general statements on the values of $\mu(G-e)$, $\muo(G-e)$, $\mud(G-e)$, and $\mut(G-e)$.  
	
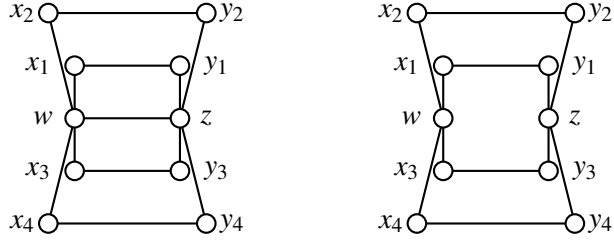
\begin{figure}[ht!]
		\begin{center}
			\begin{tikzpicture}[scale=0.7,style=thick,x=1cm,y=1cm]
				\def\vr{5pt}
				\begin{scope}[xshift=-1cm, yshift=0cm] 
					\coordinate(x1) at (0.0,1.0);
					\coordinate(x2) at (-0.5,2.0);
					\coordinate(x3) at (0.0,-1.0);
					\coordinate(x4) at (-0.5,-2.0);
					\coordinate(y1) at (2.0,1.0);
					\coordinate(y2) at (2.5,2.0);
					\coordinate(y3) at (2.0,-1.0);
					\coordinate(y4) at (2.5,-2.0);
					\coordinate(z) at (0.0,0.0);
					\coordinate(w) at (2.0,0.0);
					
					\draw (z) -- (w);  
					\draw (z) -- (x1) -- (y1) -- (w);  
					\draw (z) -- (x2) -- (y2) -- (w); 
					\draw (z) -- (x3) -- (y3) -- (w); 
					\draw (z) -- (x4) -- (y4) -- (w); 
					
					\begin{scriptsize}
						\foreach \i in {1,2,3,4} 
						{
							\draw(x\i)[fill=white] circle(\vr);
							\draw(y\i)[fill=white] circle(\vr);
						}
						\draw(z)[fill=white] circle(\vr);
						\draw(w)[fill=white] circle(\vr);
					\end{scriptsize}
					\draw (-0.6,0.0) node {$w$};
					\draw (2.5,0.0) node {$z$};
					\draw (-0.7,1.0) node {$x_1$};
					\draw (2.7,1.0) node {$y_1$};
					\draw (-1.0,2.0) node {$x_2$};
					\draw (3,2.0) node {$y_2$};
					\draw (-0.7,-1.0) node {$x_3$};
					\draw (2.7,-1.0) node {$y_3$};
					\draw (-1.0,-2.0) node {$x_4$};
					\draw (3,-2.0) node {$y_4$};
					
				\end{scope}
				\begin{scope}[xshift=6cm, yshift=0cm] 
					\coordinate(x1) at (0.0,1.0);
					\coordinate(x2) at (-0.5,2.0);
					\coordinate(x3) at (0.0,-1.0);
					\coordinate(x4) at (-0.5,-2.0);
					\coordinate(y1) at (2.0,1.0);
					\coordinate(y2) at (2.5,2.0);
					\coordinate(y3) at (2.0,-1.0);
					\coordinate(y4) at (2.5,-2.0);
					\coordinate(z) at (0.0,0.0);
					\coordinate(w) at (2.0,0.0);
					
					\draw (z) -- (x1) -- (y1) -- (w);  
					\draw (z) -- (x2) -- (y2) -- (w); 
					\draw (z) -- (x3) -- (y3) -- (w); 
					\draw (z) -- (x4) -- (y4) -- (w); 
					
					\begin{scriptsize}
						\foreach \i in {1,2,3,4} 
						{
							\draw(x\i)[fill=white] circle(\vr);
							\draw(y\i)[fill=white] circle(\vr);
						}
						\draw(z)[fill=white] circle(\vr);
						\draw(w)[fill=white] circle(\vr);
					\end{scriptsize}
					\draw (-0.6,0.0) node {$w$};
					\draw (2.5,0.0) node {$z$};
					\draw (-0.7,1.0) node {$x_1$};
					\draw (2.7,1.0) node {$y_1$};
					\draw (-1.0,2.0) node {$x_2$};
					\draw (3,2.0) node {$y_2$};
					\draw (-0.7,-1.0) node {$x_3$};
					\draw (2.7,-1.0) node {$y_3$};
					\draw (-1.0,-2.0) node {$x_4$};
					\draw (3,-2.0) node {$y_4$};
					
				\end{scope}

			\end{tikzpicture}
			\caption{Graph $H_4$ and $H_4-e$ where $e =zw$.
			}
			\label{fig:H_h}
		\end{center}
	\end{figure}

	For an integer $k$ with $k \ge 2$, let  $H_k$ be the graph defined on the vertex set $$V(H_k)= \{z, w\} \cup \{x_i, y_i: i\in [k] \}$$ such that $x_iy_izwx_i$ is a $4$-cycle for every $i \in [k]$ (see Fig.~\ref{fig:H_h}). 
	
	
	\begin{proposition} \label{prop:Hk}
		Let $e$ denote the edge $zw$ in $H_k$ for an integer $k \ge 2$.
		\begin{itemize}
			\item[(i)] $\mu(H_k)= \muo(H_k)= \mud(H_k)= \mut(H_k)= 2k$;
			\item[(ii)] $\mu(H_k-e)=k+1$ and 
			$\muo(H_k-e)=k$;
			\item[(iii)]  $\mud(H_k-e)=2$ and
			$\mut(H_k-e)=0$;  
		\end{itemize}
	\end{proposition}
	\proof (i) Let $X$ be a $\mu$-set in $H_k$. If  $X$ contains $w$ and its two neighbors $x_i$ and $x_j$, then $x_i$ and $x_j$ are not $X$-visible, a contradiction. The situation is analogous for $z$ and therefore, the condition $X\cap\{ z,w\} \neq \emptyset$ implies $|X| \leq \max\{4, k+2\} \leq   2k$. If $X\cap\{ z,w\}= \emptyset$, then $|X| \leq 2k$ automatically follows. We may infer that $\mu(H_k) \leq 2k$. On the other hand, $X= V(H_k) \setminus \{z,w\}$ is a total mutual-visibility set of cardinality $2k$ in $H_k$, and therefore, $\mut(H_k) \ge 2k$. 
	By~\eqref{eq:1}, we obtain $2k \le \mut(H_k) \le \mu(H_k)\le 2k $  and conclude $\mu(H_k)= \muo(H_k)= \mud(H_k)= \mut(H_k)= 2k$.
	\medskip
	
	(ii) Let $X$ be a mutual-visibility set in $H_k-e$. If $w \in X$, then $X$ may contain at most one vertex from $\{ x_i,y_i\}$ otherwise $w$ and $y_i$ are not $X$-visible, for every $i \in [k]$. The argument is similar for $z \in X$. Further, if both $z$ and $w$ belong to $X$, then at least one shortest path $wx_sy_sz$ must ensure the $X$-visibility of $w$ and $z$. Therefore, $X\cap\{ z,w\} \neq \emptyset$ implies $|X| \leq k+1$.
	Let us assume now $X\cap\{ z,w\} = \emptyset$. Supposing that $\{x_i,y_i, x_j, y_j\} \subseteq X$ also holds for two indices $i \neq j$, we get a contradiction as then $x_i$ and $y_j$ cannot see each other. We may infer $|X| \leq k+1 $ again. It can be checked directly that the set $Y= \{x_1, \dots , x_k, y_1\}$ is a mutual-visibility set, and we conclude that $\mu(H_k-e) = k+1$. 
	
	Now, let $X$ be an outer mutual-visibility set in $H_k-e$. If $\{z,w\} \subseteq X$, then $X=\{z,w\}$. If $w \in X$ and $z \notin X$, then $X$ may contain at most one vertex, namely $y_i$, from each pair $\{x_i, y_i\}$ and one path $wx_sy_sz$ must ensure the $X$-visibility of $w$ and $z$. Hence, $|X| \leq k$ in this case. Suppose now that $X \cap \{w,z\} = \emptyset$ and $|X|>k$. Then $X$ must contain two adjacent vertices, say $x_1$ and $y_1$, and a further vertex, say $x_2$. But then, $x_1$ and $y_2$ are not $X$-visible, which is a contradiction, and we may infer $|X| \leq k$ again. On the other hand, the $k$-element set $Y= \{x_1, \dots , x_k\}$ is an outer mutual-visibility set. This proves $\muo(H_k-e) = k$.
	\medskip
	
	(iii) If $k=2$, then $\mud (H_2-e)=\mud(C_6)= 2$ by Proposition~\ref{prop:cycles} (iii). So, we may assume $k \ge 3$. Let $X$ be a dual mutual-visibility set in $H_k-e$. Observe that the star $K_{1,k}$ induced by $N_{G-e}[w]$ is a convex subgraph in $G-e$. As $w$ has at least three neighbors, at least two of them are inside $X$, or at least two are outside $X$. In either case, these vertices have to be $X$-visible, but this is not the case if $w \in X$. This shows $w \notin X$ and, by symmetry, $z \notin X$. Suppose now that a vertex $x_i$ (resp.\ $y_i$) is in $X$. 
	Note that then $w$ and $y_i$ (resp.\ $z$ and $x_i$) are not $X$-visible, and therefore, $y_i \notin X$ (resp.\ $x_i\notin X$) would yield a contradiction. 
	Consequently, for every $i \in [k]$, it holds that $x_i \in X$ if and only if $y_i \in X$.  If $X$ contains $x_i,y_i, x_j, y_j$, then $x_i$ and $y_j$ cannot see each other. We conclude that $X$ is either empty or contains exactly two vertices, e.g., $X=\{x_1,y_1\}$. Therefore, $\mud(H_k-e)=2$.
	
	In $H_k-e$, every vertex is the middle vertex of a convex $P_3$-subgraph, and therefore, only the empty set is a total mutual-visibility set; that is, $\mut(H_k-e)=0$.  
	\qed
	
	\subsection{Mutual-visibility number under edge removal} \label{subsec:3.2}
	\begin{theorem}    
		\label{thm:edge-mu}
		If $G$ is a graph, $e \in E(G)$, and $G-e$ is  connected, then
		$$ \frac{1}{2}\, \mu(G) \leq \mu(G-e) \leq 2\,\mu(G).
		$$
		Furthermore, the lower bound is asymptotically sharp.
	\end{theorem}
	\proof Let $e=xy$ be an edge in $G$ such that $G-e$ remains connected.
	\paragraph{Upper bound.} Consider $G'= G-e$, a $\mu$-set $Y$ in $G'$, and $Y_1=Y_{G'}(x,y)$. Recall that $Y_1$ contains a vertex $z \in Y$ if $d_{G'}(z,x) \leq d_{G'}(z,y)$. 
	We obtain $G$ by adding the edge $e=xy$ to $G'$. If two vertices $u$ and $v$ are $Y$-visible in $G'$ but not in $G$, then a shortest $u,v$-path exists in $G$ that includes $e$. Then, switching the names of $u$ and $v$ if necessary, $d_{G}(u,x)<d_{G}(u,y)$ and $d_{G}(v,x)>d_{G}(v,y)$ can be observed. These relations cannot change if we refer to distances in $G'$ instead of $G$. We may conclude that if two vertices $u,v \in Y$ are not $Y$-visible in $G$, then exactly one of them belongs to $Y_1$. Therefore, any two vertices from $Y_1$ are $Y$-visible in $G$, and in turn, they are $Y_1$-visible in $G$. It proves that $Y_1$ is a mutual-visibility set in $G$. The same is true for $Y_2=Y_{G'}(y,x)$. By Observation~\ref{obs:Z-and-blocking}~(ii) we get the desired result 
	$$\mu(G) \ge \max\{|Y_1|, |Y_2|\} \ge \frac{|Y|}{2}= \frac{1}{2}\, \mu(G-e).
	$$
	\paragraph{Lower bound.} Let $X$ be a $\mu$-set in $G$ and let $X_1=X_G(x,y)$. 
	Suppose that after removing $e$, two vertices $u$ and $u'$ from $X_1$ are not $X$-visible anymore. Since $u$ and $u'$ are $X$-visible in $G$, this situation may occur only if a shortest $u,u'$-path in $G$ contains the edge $e=xy$. But then, Observation~\ref{obs:Z-and-blocking}~(i) implies that $u,u' \in X_1$ is not possible, a contradiction. Therefore, $X_1$ is a mutual-visibility set in $G-e$, and the same can be proved for $X_2= X_G(y,x)$. The lower bound readily follows if we apply Observation~\ref{obs:Z-and-blocking}~(ii) and get
	$$\mu(G-e) \ge \max\{|X_1|, |X_2|\} \ge \frac{|X|}{2}= \frac{1}{2}\, \mu(G).
	$$
	\paragraph{Sharpness.} Concerning the asymptotic sharpness of the lower bound, we show that for every $\epsilon >0$ there exists a graph $G_\epsilon$ and an edge $e_\epsilon \in E(G_\epsilon)$ such that 
	$$  \left(\frac{1}{2}+\epsilon \right)\, \mu(G_\epsilon) > \mu(G_\epsilon-e_\epsilon) .
	$$
	For a given positive $\epsilon $, let $k$ be an integer satisfying $k > \frac{1}{2\epsilon}$ and consider the graph $H_k$ and edge $e=zw$.
	Proposition~\ref{prop:Hk} (i) and (ii) then give the required inequality
	$$  \frac{\mu(H_k-e)}{\mu(H_k)}  = \frac{k+1}{2k} = \frac{1}{2} + \frac{1}{2k} < \frac{1}{2} + \epsilon.
	$$
	This finishes the proof of the theorem. \qed
	
	
	We do not have a sharp example for the upper bound in Theorem~\ref{thm:edge-mu}, but the following graph shows that the ratio $\frac{\mu(G-e)}{\mu(G)}$ can be as large as $\frac{5}{3}$. 
	Let graph $J$ be constructed from a $12$-cycle $v_1\dots v_{12}v_1$ by attaching leaves $u_1$, $u_5$, $u_9$ to the vertices $v_1$, $v_5$, $v_9$, respectively. Let $X$ be a $\mu$-set in $J$. Since the $12$-cycle is a convex subgraph, $X$ contains at most three vertices from it. If a leaf, say $u_1$, is in $X$, then $X$ contains at most two vertices from the cycle. Indeed, if $\{v_i, v_j, v_k\} \subseteq X$ and $1 \leq  i<j<k$, then $v_j$ cannot see $u_1$. It can be proved similarly that in case of $\{u_1, u_5\} \subseteq X$, the $\mu$-set contains at most one vertex from the cycle; and if $\{u_1, u_5, u_9\} \subseteq X$, then no vertex from the cycle belongs to $X$. Furthermore, $X=\{u_1, u_5, u_9\}$ is clearly a mutual-visibility set in $J$. This proves $\mu(J)=3$. By removing the edge $e=v_2v_3$ from $J$, a tree with five leaves is obtained. Therefore, $\mu(J-e)=5=\frac{5}{3}\, \mu(J)$. Similar constructions with a ratio of $5/3$ can be obtained from longer cycles as well.

	\subsection{Outer mutual-visibility under edge removal} 
	\begin{theorem}    
		\label{thm:edge-outer}
		If $G$ is a graph, $e \in E(G)$, and $G-e$ is connected, then
		$$ \frac{1}{6}\, \muo(G) \leq  \muo(G-e) \leq 2\,\muo(G)+1.
		$$
	\end{theorem}
	\proof 
	Let $G$ be a graph and $e=xy$ an edge in $G$ such that $G-e$ is connected. 
	\paragraph{Upper bound.} Choose an arbitrary $\muo$-set $Y$ in $G'=G-e$.
	In the proof, we use an approach similar to that in the proof of Theorem~\ref{thm:edge-mu} and define $Y_1=Y_{G'}(x,y)$.
	We consider two cases.
	\medskip
	
	{\it Case 1.} Assume first that $x,y \in Y$ and prove that $Y'=Y \setminus \{x,y\}$ is an outer mutual-visibility set in $G=G'+e$. Let $u$ be a vertex from $ Y' $ and $v$ be an arbitrary vertex in $G$. If no shortest $u,v$-path in $G$ contains the edge $e$, then $u$ and $v$ are $Y$-visible in $G$ via the same path as in $G'$. 
	If a shortest $u,v$-path in $G$ contains the edge $e$ and $d_{G'} (u,x) < d_{G'}(u,y)$, then the distance of $u$ and $v$ in $G$ is $d_G(u,v)= d_G(u,x)+1+d_G(y,v)$. 
	Since $Y$ is an outer mutual-visibility set in $G'$ and $x,y \in Y$, there exist a shortest $u,x$-path $P^1$ and a shortest $y,v$-path $P^2$ in $G'$ which contain no internal vertices from $Y$. Observe that $P^1$ and $P^2$ remain isometric paths in $G$. Then $P^1$, $e$, and $P^2$ together form a shortest $u,v$-path in $G$, which does not contain internal vertices from $Y'$. If $d_{G'} (u,x) > d_{G'}(u,y)$, the proof is similar; and if $d_{G'} (u,x) = d_{G'}(u,y)$, these distances remain the same in $G$ and a shortest $u,v$-path in $G$ cannot contain $e$. Hence, $Y'$ is always an outer mutual-visibility set in $G$ and therefore, $\muo(G) \ge \muo(G-e)-2$. Since $\muo(G) \ge 2$, we can conclude $\muo(G-e) \leq \muo(G)+2 < 2 \muo(G)+1$. 
	\medskip
	
	{\it Case 2.} We now assume that $|Y\cap \{x,y\}| \leq 1$ and show that $Y_1'=Y_1 \setminus \{x\}$ is an outer mutual-visibility set in $G$. By symmetry, the same holds for $Y_2'=Y_{G'}(y,x)\setminus \{y\}$.
	
	Let $u \in Y_1'$ and $v$ be an arbitrary vertex from $G$. Suppose, to the contrary, that $u$ and $v$ are not $Y_1'$-visible in $G$. As these vertices are $Y_1'$-visible in $G'$, their distance is smaller in $G$ than in $G'$. Consequently, every shortest $u,v$-path in $G$ contains $e$. Since $u \in Y$, there exists a shortest $u,x$-path $P^1$ in $G'$ such that no internal vertex of $P^1$ belongs to $Y$. Consequently, no internal vertex of $P^1$ is from $Y_1'$. Let $P^2$ be a shortest $y,v$-path in $G'$. Since each shortest $u,v$-path contains $e$, every vertex from $P^2$ is (strictly) nearer to $y$ than to $x$. Thus $V(P^2) \cap Y_1'= \emptyset$. As $x \notin Y_1'$ also holds, the path obtained from $P^1$, $e$, and $P^2$ is a shortest $u,v$-path in $G$ and no internal vertex of it belongs to $Y_1'$. We may infer that $Y_1'$ is an outer mutual-visibility set in $G$, and the same holds for $Y_2'$. Our assumption $|Y\cap \{x,y\}| \leq 1$ implies $|Y_1'| + |Y_2'| \ge |Y|-1$ and hence we get the desired relation
	$$\muo(G) \ge \max\{|Y_1'|, |Y_2'|\} \ge \frac{|Y|-1}{2}= \frac{\mu(G-e)-1}{2}. 
	$$
	
	\paragraph{Lower bound.} Suppose now that $X$ is a $\muo$-set in $G$. We will show that one can always find a subset of $X$ of size at least $|X|/6$ that is an outer mutual-visibility set in $G'=G-e$. Define first $X_1=X_{G}(x,y)$ and $X_2=X_{G}(y,x)$. Then $|X_1| + |X_2| \ge |X|$ and, by symmetry, we may suppose $|X_1| \ge |X|/2$. We partition $X_1$ into two sets and then, prove a series of claims. Let
	\begin{align*}
		X_{1,1} & =\{v \in X_1: d_{G-e}(v,y) \leq d_{G-e}(v,x)+1\} \\
		X_{1,2} & =\{v \in X_1: d_{G-e}(v,y) \geq d_{G-e}(v,x)+2\}. 
	\end{align*}
	\noindent{\it Claim A.}\enskip  If $u \in X_1$ and $v$ is a vertex in $G$ with $d_G(v,x) \leq d_G(v,y)$, then $u$ and $v$ are $X_1$-visible in $G-e$.\\
	\noindent{\it Proof.}\enskip By Lemma~\ref{lem:subset-closed} $X_1$ is an outer mutual-visibility set in $G$, and then $u$, $v$ are $X_1$-visible in $G$. Thus, a shortest $u,v$-path $P$ exists in $G$ that contains no internal vertex from $X_1$. If edge $e=xy$ lies on $P$, then, since $P$ is an isometric subgraph in $G$, one of $u$ and $v$ is (strictly) closer to $y$ than to $x$, which contradicts the conditions in the claim. Therefore, $e$ does not belong to $P$, and the two vertices remain $X_1$-visible via $P$ in $G-e$. \smallqed
	\medskip
	
	\noindent{\it Claim B.}\enskip   $X_{1,2}$ is an outer mutual-visibility set in $G-e$.\\
	\noindent{\it Proof.}\enskip Suppose, to the contrary, that two vertices $u \in X_{1,2}$ and $v \in V(G)$ are not $X_{1,2}$-visible in $G-e$. By Claim~A, $v \notin X_{1,2}$. Let $P$ be a shortest $u,v$-path in $G-e$ such that $|V(P) \cap X_{1,2}|$ is minimum. By Claim~A, this path $P$ cannot contain more than one internal vertex from $X_{1,2}$. Hence, $V(P) \cap X_{1,2} =\{u, u_1\}$ for  a vertex $u_1 \in X_{1,2}$. Let $u_2$ be the neighbor of $u_1$ that follows it (towards $v$)  on the $u,v$-path $P$. By our assumption, $u$ and $v$ are not $X_{1,2}$-visible and thus, $u$ and $u_2$  are not $X_{1,2}$-visible in $G-e$, either. Claim~A then implies $d_G(u_2,x) > d_G(u_2,y)$ and, in turn, $d_{G-e}(u_2,x) > d_{G-e}(u_2,y)$ follows. 
	Also, the removal of $e$ does not affect the distance between $u_2$ and $y$; that is, $d_{G-e}(u_2,y)=d_G(u_2,y)$. As there is an $x,u_2$-path in $G-e$ through $u_1$, we have $d_{G-e}(u_2,x) \leq d_{G-e}(u_1,x) +1 $. Similarly, $d_{G-e}(u_1,y)$ can be estimated as $d_{G-e}(u_1,y) \leq d_{G-e}(u_2,y) +1 $. This yields
	$$ d_{G-e}(u_1,y) \leq d_{G-e}(u_2,y) +1 <d_{G-e}(u_2,x)+1 \leq d_{G-e}(u_1,x) +2 
	$$
	that contradicts the assumption $u_1 \in X_{1,2}$. Thus $X_{1,2}$ is an outer mutual-visibility set in $G-e$. \smallqed
	\medskip
	
	Before the last claim is stated, we specify a subset $X_{1,1}'$ of $X_{1,1}$. By Claim~A, $X_{1,1}$ is a mutual-visibility set in $G-e$. We will remove blocking vertices from $X_{1,1}$ as defined in Section~\ref{sec:pre}. First, we choose a shortest $y,u$-path $P_{y,u}$, for each $u \in X_{1,1}$, that contains the least internal vertices from $X_{1,1}$. Then, we remove a vertex $u'$ from $X_{1,1}$, if $u'$ is an internal vertex of $P_{y,u}$ for some $u \in X_{1,1}$. That is, we consider $X_{1,1}'=X_{1,1} \setminus B_{X_{1,1},y}$. Observation~\ref{obs:Z-and-blocking} (iv) shows that $|X_{1,1}'| \ge |X_{1,1}|/2$.
	\medskip
	
	\noindent{\it Claim C.}\enskip   $X_{1,1}'$ is an outer mutual-visibility set in $G-e$.\\
	\noindent{\it Proof.}\enskip Suppose, to the contrary, that two vertices $u \in X_{1,1}'$ and $v \in V(G)$ are not $X_{1,1}'$-visible in $G-e$.  Then, there exists a shortest $u,v$-path $P$ in $G-e$ that contains exactly one internal vertex, say $u_1$ from $X_{1,1}'$. Let $u_2$ be the neighbor of $u_1$ that follows $u_1$ on the path $P$ towards $v$. Clearly, $u$ and $u_2$ are not $X_{1,1}'$-visible in $G-e$, either, and Claim A implies $d_G(u_2,y) < d_G(u_2,x)$.
	
	Since $u$ and $u_2$ are $X_{1,1}'$-visible in $G$ but not in $G-e$, a shortest $u,u_2$-path in $G$ contains $e$. Then the two distances can be compared as
	\begin{equation} \label{eq:2}
		d_{G-e} (u,u_2) \ge d_{G} (u,u_2) = d_G(u,x) + 1+ d_{G} (u_2,y)=   d_{G-e}(u,x) + 1+ d_{G-e} (u_2,y),  
	\end{equation}
	where we also used the fact that removing $e=xy$ cannot change $d_G(u,x)$ and $d_G(u_2,y)$ as  $d_G(u,x) \leq  d_G(u,y)$ and $d_G(u_2,y) < d_G(u_2,x)$. Since $u \in X_{1,1}$ also holds, it is true that
	\begin{equation} \label{eq:3}
		d_{G-e}(u,x) +1 \ge d_{G-e} (u,y).
	\end{equation}
	From (\ref{eq:2}) and (\ref{eq:3}), we obtain
	\begin{equation} \label{eq:4}
		d_{G-e}(u,u_2) \ge  d_{G-e} (u,y) + d_{G-e} (u_2,y).
	\end{equation}
	Since there is a $u,u_2$-path through $y$, inequality~(\ref{eq:4}) holds with equality and proves that there is a shortest $u,u_2$-path through $y$  in $G-e$. By Observation~\ref{obs:Z-and-blocking} (iv), $u$ and $y$ became $X_{1,1}'$-visible after the removal of the blocking vertices. Since $u$ and $u_2$ are not $X_1$-visible, Claim A implies $d_{G-e}(u_2,x) > d_{G-e}(u_2, y)$ and therefore, $y$ and $u_2$ are  $X_{1,1}'$-visible in $G-e$. We may infer the $X_{1,1}'$-visibility of $u$ and $u_2$, a contradiction that verifies the claim. \smallqed
	
	To establish the lower bound, we first recall the assumption $|X_1| \ge |X|/2$. Claims B and C show that both $X_{1,2}$ and $X_{1,1}'$ are outer mutual-visibility sets in $G-e$. Observation~\ref{obs:Z-and-blocking} (ii) gives $|X_{1,1}'| \ge |X_{1,1}|/2$. Therefore, we get
	$$  \muo(G-e) \ge \max\left\{|X_{1,2}|, \frac{|X_{1,1}|}{2}\right\} \ge \frac{|X_{1}|}{3} \ge \frac{|X|}{6} = \frac{1}{6}\, \muo(G) 
	$$   
	as stated in the theorem.
	\qed 
	\subsection{Total and dual mutual-visibility under edge removal}
	\begin{theorem}
		\label{thm:edge-mut}
		Suppose that $G$ is a graph, $e$ is an edge in $G$, and $G-e$ is connected. Then the following statements are valid.
		\begin{itemize}
			\item[(i)] $  \mut(G-e) \leq \mut(G) +2$; and this 
			inequality is sharp.
			\item[(ii)] There are no constants $a>0$ and $b$  such that $a\cdot \mut(G) +b \leq \mut(G-e)$ holds for every such graph $G$ and edge $e \in E(G)$.
		\end{itemize}
	\end{theorem}
	\proof
	(i) First, consider the graph $G'=G-e$ where $e=xy$, and let $Y$ be a $\mut$-set in $G'$. We want to prove that $Y'=Y \setminus \{x,y\}$ is a total mutual-visibility set in $G=G'+e$. By Lemma~\ref{lem:distance two}, it suffices to show that every $u,v\in V(G)$ with $d_G(u,v)=2$ are $Y'$-visible. Let $P$ be a shortest $u,v$-path in $G$ where $d_G(u,v)=2$. If $P$ does not go through $e$, then $d_{G'}(u,v)=d_{G}(u,v)= 2$ and since $u,v$ are $Y'$-visible in $G'$, they remain $Y'$-visible in $G$. If $P$ contains the edge $e$, then the middle vertex of $P$ is either $x$ or $y$. So the middle vertex of $P$ is not from $Y'$, and $u$ and $v$ are $Y'$-visible in $G$. Therefore, $Y'$ is a total mutual-visibility set in $G$ and  $\mut(G) \ge \mut(G-e) -2$. 
	
	The sharpness is shown by all cycles $C_n$ with $n \ge 5$ as $\mut(C_n)=0$ and $\mut(C_n-e)=\mut(P_n)=2$ for every edge $e$ of $C_n$. To obtain more general examples for the sharpness, we may start with a cycle $C_n \colon v_1v_2\dots v_nv_1$, with $n \ge 5$, and attach some leaves to some of the vertices $v_3, \dots, v_n$. Let $\ell$ be the total number of leaves attached and denote by $C_n^\ell$ the graph constructed this way. Let $X$ be a $\mut$-set in $C_n^\ell$. As the $n$-cycle is a convex subgraph in $C_n^\ell$, Lemma~\ref{lem:convex} implies that $X$ contains no vertices from the cycle. Thus, $\mut(C_n^\ell) \leq \ell$. On the other hand, the set of leaves is a total mutual-visibility set in $C_n^\ell$ and therefore, $\mut(C_n^\ell)= \ell$. 
	The removal of the edge $e=v_1v_2$ results in a tree with $\ell+2$ leaves and so $\mut(C_n^\ell-e) =\ell+2$.
	\medskip
	
	(ii) Proposition~\ref{prop:Hk} (i) and (iii) show that 
	$\mut(H_k)=2k$ and $\mut(H_k-e)=0$ for every $k \ge 2$. For every two constants $a$, $b$, we may choose an integer $k >\max\{\frac{-b}{2a},1\}$ and then $a\cdot \mut(H_k) +b > \mut(H_k-e)$ follows. This verifies $(ii)$.  \qed
	\begin{proposition}
		\label{prop:edge-mud}
		Suppose that $G$ is a graph, $e$ is an edge in $G$, and $G-e$ is connected. 
		\begin{itemize}
			\item[(i)]  There are no constants $a>0$ and $b$ such that $a\cdot \mud(G) +b \leq \mud(G-e)$ holds for every $G$ and $e$.
			\item[(ii)] There are no constants $a>0$ and $b$ such that $a\cdot \mud(G) +b \geq \mud(G-e)$ holds for every $G$ and $e$.
		\end{itemize}
	\end{proposition}
	\proof
	(i) The non-existence of the lower bound follows from Proposition~\ref{prop:Hk}~(i) and (iii) if we choose an integer $k >\max\{ \frac{2-b}{2a},1\}$ and consider $H_k$.
	
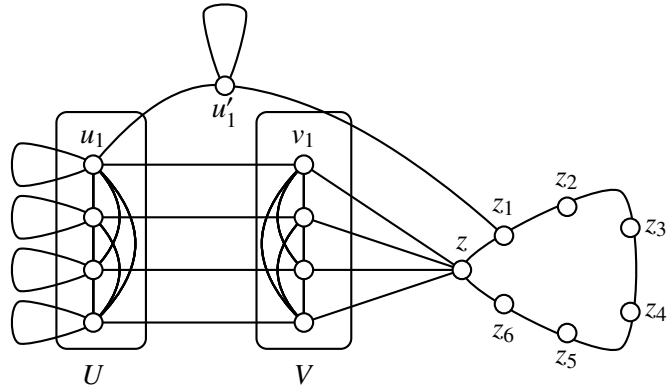
\begin{figure}[ht!] 
		\centering
		\begin{tikzpicture}
			[scale=0.7,style=thick,x=1cm,y=1cm]
			\def\vr{5pt}
			\begin{scope}[xshift=-1cm, yshift=0cm] 
				\coordinate(u1) at (0.0,0.0);
				\coordinate(u2) at (0.0,-1.0);
				\coordinate(u3) at (0.0,-2.0);
				\coordinate(u4) at (0.0,-3.0);
				\coordinate(v1) at (4.0,0.0);
				\coordinate(v2) at (4.0,-1.0);
				\coordinate(v3) at (4.0,-2.0);
				\coordinate(v4) at (4.0,-3.0);
				\coordinate(z) at (7,-2);
				\coordinate(z1) at (7.8,-1.35);
				\coordinate(z6) at (7.8,-2.65);
				\coordinate(z2) at (9,-0.8);
				\coordinate(z5) at (9,-3.2);
				\coordinate(z3) at (10.2,-2.8);
				\coordinate(z4) at (10.2,-1.2);
				\coordinate(u1') at (2.5,1.5);
				
				\draw plot [smooth, tension=1] coordinates {(u1) (1.3,1.2) (u1')};
				\draw plot [smooth, tension=1] coordinates {(u1') (5,0.7) (z1)};
				\draw plot [smooth cycle] coordinates {(0,0) (-1.4,0.4) (-1.4,-0.4)};
				\draw plot [smooth cycle] coordinates {(0,-1) (-1.4,-0.6) (-1.4,-1.4)};
				\draw plot [smooth cycle] coordinates {(0,-2) (-1.4,-1.6) (-1.4,-2.4)};
				\draw plot [smooth cycle] coordinates {(0,-3) (-1.4,-2.6) (-1.4,-3.4)};
				\draw plot [smooth cycle] coordinates {(7,-2) (10,-3.5) (10,-0.5)};
				\draw plot [smooth cycle] coordinates {(2.5,1.5) (2.1,2.9) (2.9,2.9)};
				\foreach \i in {1,2,3,4} 
				{
					\draw(z) -- (v\i);
					\draw(v\i) -- (u\i);
					\draw(v1) -- (v2) -- (v3) -- (v4);
					\draw(u1) -- (u2) -- (u3) -- (u4);
					\draw plot [smooth, tension=1] coordinates {(u1) (0.8,-1.5) (u4)};
					\draw plot [smooth, tension=1] coordinates {(u1) (0.5,-1) (u3)};
					\draw plot [smooth, tension=1] coordinates {(u2) (0.5,-2) (u4)};
					\draw plot [smooth, tension=1] coordinates {(v1) (3.2,-1.5) (v4)};
					\draw plot [smooth, tension=1] coordinates {(v1) (3.5,-1) (v3)};
					\draw plot [smooth, tension=1] coordinates {(v2) (3.5,-2) (v4)};
				}

				\begin{scriptsize}
					\foreach \i in {1,2,3,4} 
					{
						\draw(u\i)[fill=white] circle(\vr);
						\draw(v\i)[fill=white] circle(\vr);
					}
					\foreach \i in {1,2,3,4,5,6} 
					{
						\draw(z\i)[fill=white] circle(\vr);
					}
					\draw(z)[fill=white] circle(\vr);
					\draw(u1')[fill=white] circle(\vr);
				\end{scriptsize}
				
				\draw (7, -1.5) node {$z$};
				\draw (7.8, -0.8) node {$z_1$};
				\draw (9, -0.3) node {$z_2$};
				\draw (10.7, -1.2) node {$z_3$};
				\draw (10.7, -2.8) node {$z_4$};
				\draw (9, -3.7) node {$z_5$};
				\draw (7.8, -3.2) node {$z_6$};
				\draw (2.5,1) node {$u_1'$};
				\draw (0,0.5) node {$u_1$};
				\draw (4,0.5) node {$v_1$};
			\end{scope}
			\draw[rounded corners] (0, -3.5) rectangle (-1.7, 1);
			\draw (-1,-4) node {$U$};
			\draw[rounded corners] (2.1, -3.5) rectangle (3.9, 1);
			\draw (3.,-4) node {$V$};
			
		\end{tikzpicture}
		\caption{A schematic drawing of graph $L_4$, where ovals represent $7$-cycles.}
		\label{fig:L_4}
	\end{figure}

	(ii) To prove the non-existence of the upper bound, we construct a graph $L_k$ for every $k \ge 2$ as follows. Let $U=\{u_1, \dots, u_k\}$ and $V=\{v_1, \dots, v_k\} $ and take a further vertex $z$. For the edges, take a matching $\{u_iv_i: i \in [k]\}$, add all edges inside $U$ and inside $V$, and make $z$ adjacent to all vertices in $V$. Then we put vertex disjoint $7$-cycles onto all vertices in $U \cup \{z\}$. That is, for every  $w \in U \cup \{z\}$, we take six new vertices which form a $C_7$ with $w$. Let this $7$-cycle be $zz_1\dots z_6z$ for the vertex $z$. Finally, we add a new vertex $u_1'$, two incident edges $u_1u_1'$ and $z_1u_1'$, and put a $7$-cycle onto $u_1'$. See Fig.~\ref{fig:L_4} for an illustration.
	
	We first show that $\mud(L_k)=0$ holds for every $k$. Let $X$ be a $\mud$-set of $L_k$. Every induced $7$-cycle $C$ is a convex subgraph in $L_k$. Lemma~\ref{lem:convex} implies that $V(C) \cap X$ is a dual mutual visibility set in $C$. Since $\mud(C_7)=0$, no vertex from an induced $7$-cycle belongs to $X$ in $L_k$. In particular, $X \subseteq V$. On the other hand, every vertex $v_i \in V$ lies on the unique shortest $zu_i$-path and therefore, $u_i \notin X$ and $z \notin X$ imply $v_i \notin X$. This shows $X= \emptyset$ and $\mud(L_k)= 0$.
	
	Now, let $e=zz_6 $ and consider the graph $L_k-e$. It is straightforward to check that $Y=V \cup \{z\}$ is a $(k+1)$-element dual mutual-visibility set in $L_k-e$. Remark that vertices $u_i$ and $z_j$ can see each other in $L_k-e$ via a shortest path that goes through $u_1'$. Then $\mud(L_k-e)\ge k+1$. (Also, it is not hard to see that $\mud(L_k-e)=k+1$.) Thus $a\cdot \mud(L_k) +b < \mud(L_k-e)$ holds when  $k \ge \max\{b, 2\}$ and we are done.
	\qed
	
	\begin{remark}
		Since graphs $H_k$ in  Proposition~\ref{prop:Hk} are bipartite graphs of diameter $3$, the statements in Theorem~\ref{thm:edge-mut} (ii) and Proposition~\ref{prop:edge-mud} (i) remain true over the class of bipartite graphs of diameter $3$.
	\end{remark}
	
	\section{Realizability of visibility invariants in terms of the order} \label{sec:realizability}
	
	Our goal in this section is to characterize $(\sigma, n)$-realizable pairs $(p,q)$ for every visibility invariant $\sigma \in \{\mu, \muo, \mud, \mut\}$. Recall that $(p,q)$ is $(\sigma, n)$-realizable if there exists a graph $G$ with $\sigma(G)=p$ and $n(G)=q$. These results are applied in the next section when visibility invariants under vertex removal are studied.
	
	\begin{proposition}
		\label{prop:realizability}
		Let $p$ and $q$ be two integers satisfying $0 \leq p \leq q$ and $1 \leq q$.
		\begin{itemize}
			\item[(i)] The following statements are equivalent:
			\begin{itemize}
				\item[$\circ$] $(p,q)$ is $(\mu, n)$-realizable;
				\item[$\circ$] $(p,q)$ is $(\muo, n)$-realizable;
				\item[$\circ$] $p \ge 2$ or $(p,q)=(1,1)$.                     
			\end{itemize}
			\item[(ii)] The pair $(p,q)$ is $(\mut, n)$-realizable if and only if 
			\begin{itemize}
				\item[$\circ$] $p \ge 2$; or
				\item[$\circ$] $p=1$ and $q\ge 6$; or $(p,q)=(1,1)$; or
				\item[$\circ$] $p=0$ and $q\ge 5$.
			\end{itemize}
			\item[(iii)] The pair $(p,q)$ is $(\mud, n)$-realizable if and only if 
			\begin{itemize}
				\item[$\circ$] $p \ge 2$; or
				\item[$\circ$] $p=1$ and $q\ge 8$; or $(p,q)=(1,1)$; or
				\item[$\circ$] $p=0$ and $q\ge 7$.     
			\end{itemize}
		\end{itemize}
	\end{proposition}
	\proof First recall that if $G$ is a tree, the value of $\sigma(G)$ equals the number of leaves in $G$ for each $\sigma \in \{\mu, \muo, \mud, \mut\}$ \cite{distefano-2022, CiDiDrHeKlYe-2023}. Therefore, $(p,q)$ is always $(\sigma,n)$-realizable by a tree when $2 \leq p < q$. Also, the cases with $1 \leq p=q$ can be easily settled as $\sigma(K_q)=q$ for every positive integer $q$. Clearly, $(0,1)$ is not $(\sigma, n)$-realizable. What remains to consider are the cases with $p \in \{0,1\}$ and $q \ge 2$.
	
	(i) In a graph $G$ of order $q \ge 2$, any two vertices $u,v$ form a mutual-visibility set, and if the distance of $u$ and $v$ equals the diameter of the graph, then $\{u,v\}$ is also an outer mutual-visibility set. Therefore, $\mu(G) \ge 2$ and $\muo(G) \ge 2$ always hold. In turn, if $q\ge 2$, no pair $(0,q)$ or $(1,q)$ is $(\sigma, n)$-realizable for $\sigma \in \{\mu, \muo\}$. This finishes the proof for (i).
	
	(ii) By Proposition~\ref{prop:cycles} (iv), $\mut(C_s) =0$ for every $s \ge 5$. Thus $(0,q)$ is $(\mut,n)$-realizable when $q \ge 5$. For pairs $(1,q)$ with $q \ge 6$, we take a graph $C_{q-1}^+$ obtained from $C_{q-1}$ by attaching a leaf to it. As $C_{q-1}$ is a convex subgraph in $C_{q-1}^+$, and $\mut(C_{q-1})=0$, a $\mut$-set $X$ of $C_{q-1}^+$ cannot contain any vertex from $C_{q-1}$. On the other hand, the leaf of $C_{q-1}^+$ forms a total mutual-visibility set in $C_{q-1}^+$. It proves that $(1,q)$ is $(\mut, n)$-realizable for all $q \ge 6$. 
	
	To finish the proof for (ii), it suffices to show that $(0,q)$ is not $(\mut,n)$-realizable if $q \in \{2,3,4\}$ and $(1,q)$ is not $(\mut,n)$-realizable if $q \in \{2,3,4,5\}$. The cases with $q \in \{2,3,4\}$ are easy to check directly. Consider then the pair $(1,5)$. Let $G$ be a connected graph of order $5$. If $G$ has at least two leaves, they form a total mutual-visibility set, and $\mut(G) \ge 2$ follows.  
	
	Now, assume that $G$ has at most one leaf.
	If $\Delta(G) =4$, then $G$ has a universal vertex and $\mut(G) \ge 4$. 
	If $\Delta(G) =3$, let $V(G)=\{v, v_1, v_2, v_3, u\}$ and  $N_G(v)=\{v_1, v_2, v_3\}$. Assume that $\deg(u)=1$, say $N(u)=\{ v_1\}$. Then $\{v_2, v_3, u\}$ is a total mutual-visibility set and $\mut(G) \ge 3.$ In the case of $\deg(u) \ge 2$, there are two neighbors of $v$ which form a total mutual-visibility set in $G$ and $\mut(G) \ge 2$.
	Finally, if $\Delta(G) =2$, then $G$ is a $5$-cycle and $\mut(G)=0$.
	We may conclude that $(1,5)$ is not $(\mut,n)$-realizable and that $C_5$ is the only graph that realizes $(0,5)$. 
	
	(iii)
	By Proposition \ref{prop:cycles} (iii), we have $\mud(C_s) =0$ for every $s \ge 7$. Hence $(0,q)$ is $(\mud,n)$-realizable when $q \ge 7$. 
	For pairs $(1,q)$ with $q \ge 8$, consider again the graph $C_{q-1}^+$. 
	Let $X$ be a $\mud$-set in $C_{q-1}^+$. As $C_{q-1}$ is a convex subgraph in $C_{q-1}^+$ and $\mud(C_{q-1})=0$, Lemma~\ref{lem:convex} implies $X \cap V(C_{q-1})= \emptyset$. Since the leaf in $C_{q-1}^+$ forms a dual mutual-visibility set, we may infer
	$\mud(C_{q-1}^+)= 1$. It yields that $(1,q)$ is $(\mud,n)$-realizable when $q \ge 8$. 
	
	In the proof of part (ii), we have seen that $(0,q)$ and $(1,q)$ are not $(\mut,n)$-realizable if $q \in \{2,3,4\}$.  
	Moreover, $(0,5)$ is $(\mut,n)$-realizable by only $C_5$ and $(1,5)$ is not $(\mut,n)$-realizable. As $\mut(G) \le \mud(G)$ for every  graph $G$ and $\mud(C_5)=2$, we deduce that
	$(0,q)$ and $(1,q)$ are not $(\mud,n)$-realizable when $q \in \{2,3,4,5\}$.
	
	We now show that $(0,6)$ and $(1,6)$ are not $(\mud,n)$-realizable. Let $G$ be a connected graph of order~$6$.  
	\begin{itemize}
		\item[(a)] If $G$ has at least two leaves, then the set of its leaves is a dual mutual-visibility set. Then $\mud(G) \ge 2$.
		\item[(b)] If $\Delta(G) =5$, then $G$ has a universal vertex $v$ and, as $N(v)$ is a dual mutual-visibility set,  $\mud(G) \ge 5$. 
		\item[(c)] If $\Delta(G) \le 4$ and there is a partition of $V(G)$ into $A$ and $B$ such that $G[A]$ is either $K_2$ or $K_3$, and $G[B]$ is a graph of diameter at most two, then by Lemma \ref{lem:lower}, we have $\mud(G) \ge 2$.
	\end{itemize}
	By checking all $112$ connected graphs on six vertices (see e.g.\ \cite{cvetkovic}), it can be seen that all but one of them satisfy at least one property from (a)-(c). The exception is $C_6$ and since $\mud(C_6)=2$, we may conclude that $(0,6)$ and $(1,6)$ are not $(\mud, n)$-realizable.

	It remains to prove that $(1,7)$ is not $(\mud,n)$-realizable. 
	Suppose, to the contrary, that $G$ is a graph of order $7$ with $\mud(G)=1$ and that $\{x\}$ is a $\mud$-set in $G$. It follows from Lemma~\ref{lem:mud-1} that $\mut(G-x)=0$ and, equivalently, every vertex is the middle vertex of a convex path $P_3$ in $G-x$. Checking again those $112$ connected graphs on six vertices, we find that only the $6$-cycle satisfies the condition. Hence,  $G-x$ is a $6$-cycle. Let $v_1, \dots, v_6$ be the vertices of $G-x$ (in natural order). If there exist two consecutive vertices $v_i$ and $v_{i+1}$ on the cycle such that both are adjacent to $x$, then $\{x,v_i, v_{i+1}\}$ is a dual mutual-visibility set in $G$ that contradicts $\mud(G)=1$. If there exist two consecutive vertices $v_i$ and $v_{i+1}$ on the cycle such that both are non-adjacent to $x$, then $V(G) \setminus \{v_i, v_{i+1}\}$ is an isometric subgraph of $G$ and Lemma~\ref{lem:lower} implies $\mud(G) \ge 2$, a contradiction. Hence, for every $i \in [6]$, exactly one of the vertices $v_i$ and $v_{i+1}$ is adjacent to $x$. Without loss of generality, we may assume that $x$ is adjacent to $v_1$, $v_3$, and $v_5$ in $G$. 
	It is easy to check that $v_1,v_2,v_3$ form a dual mutual-visibility set and $\mud(G) \ge 3$, a contradiction again. We conclude that $\mud(G)=1$ is not possible if $G$ is a graph of order $7$.  
	This finishes the proof of the proposition.
	\qed
	\section{Vertex removal} \label{sec:vertex-removal}
	This section is divided into two parts. First, we show that the visibility invariants can decrease arbitrarily by the removal of a vertex $x$. In the second part, we prove that $\mu(G-x)$ is bounded by $2\mu(G)$ from above but there are no similar upper bounds for $\mud(G-x)$ and $\mut(G-x)$.
	\subsection{Decreasing visibility invariants by vertex removal}
	
	We point out that the decrease in the visibility invariants may be arbitrary when a vertex is removed from a graph.
	
	\begin{proposition}
		\label{prop:decrease_vertex}
		If\/ $\sigma \in\{\mu, \muo, \mud, \mut\}$, and $p,q$ are integers satisfying $2 \leq p \leq q$, then there exists a graph $G$ and a vertex $x \in V(G)$ such that  $\sigma(G)=q$ and $\sigma(G-x)=p$.  
	\end{proposition}
	\proof If $2 \leq p=q$, the statement is verified by the example of $S(K_{1,p})$  that is obtained from a star $K_{1,p}$ by subdividing each edge once. Since $S(K_{1,p})$ is a tree with $p$ leaves,  $\sigma(S(K_{1,p}))=p$ for each $\sigma \in \{\mu, \muo, \mud, \mut\}$. Removing a leaf $x$ from $S(K_{1,p})$, we obtain a tree with $p $ leaves. Then $\sigma(S(K_{1,p})-x)=p$ as required.
	
	By Proposition~\ref{prop:realizability}, for every pair $(p,q)$ with $2 \leq p < q$, and for each $\sigma \in \{\mu, \muo, \mud, \mut\}$, we can choose a non-complete graph $G'$ with  $\sigma(G')= p$ and $n(G')=q$. Construct $G$ from $G'$ by adding a universal vertex $x$. Since $G$ is not a complete graph, we have $\mu(G) \leq n(G)-1=q$. On the other hand, $V(G) \setminus \{x\}$ is a total mutual-visibility set in $G$ and then $\mut(G) \ge q$. The inequalities in \eqref{eq:1} then imply $\mut(G)=\mud(G)=\muo(G)=\mu(G)=q$. In particular, $\sigma(G)=q$. By construction, $G'=G-x$ and  $\sigma(G-x)=p$.
	This completes the proof.
	\qed
	
	Proposition~\ref{prop:decrease_vertex} immediately implies the non-existence of the considered lower bounds when a vertex is removed from a graph.
	\begin{corollary}
		\label{cor:lower-vertex}
		The following statement holds for each visibility invariant $\sigma \in \{\mu, \muo, \mud, \mut\}$. There are no constants $a>0$ and $b$ such that $a\cdot \sigma(G) +b \leq \sigma(G-x)$ for every $G$ and $x \in V(G)$.  
	\end{corollary} 
	
	\subsection{Increasing visibility invariants by vertex removal}
	In this subsection, we prove the general upper bound  $\mu(G-x) \leq 2\mu(G)$. A further result shows that   $\mud(G-e)$ and $\mut(G-e)$ can be arbitrarily large compared to $\mud(G)$ and $\mut(G)$, respectively. We first introduce an infinite family of graphs.
	
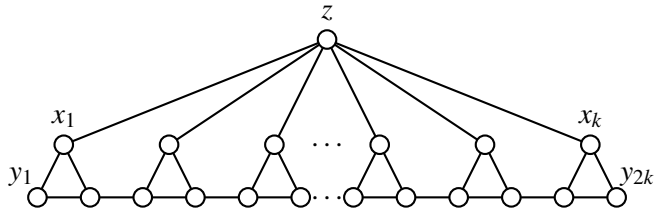
\begin{figure}[ht!]
		\begin{center}
			\begin{tikzpicture}[scale=0.7,style=thick,x=1cm,y=1cm]
				\def\vr{5pt}
				\begin{scope}[xshift=-1cm, yshift=0cm] 
					\coordinate(z) at (0.5,5.0);
					\coordinate(x1) at (-4.5,3.0);
					\coordinate(x2) at (-2.5,3.0);
					\coordinate(x3) at (-0.5,3.0);
					\coordinate(x4) at (1.5,3.0);
					\coordinate(x5) at (3.5,3.0);
					\coordinate(x6) at (5.5,3.0);
					\coordinate(y1) at (-5.0,2.0);
					\coordinate(y2) at (-4.0,2.0);
					\coordinate(y3) at (-3.0,2.0);
					\coordinate(y4) at (-2.0,2.0);
					\coordinate(y5) at (-1.0,2.0);
					\coordinate(y6) at (0.0,2.0);
					\coordinate(y7) at (1.0,2.0);
					\coordinate(y8) at (2.0,2.0);
					\coordinate(y9) at (3.0,2.0);
					\coordinate(y10) at (4.0,2.0);
					\coordinate(y11) at (5.0,2.0);
					\coordinate(y12) at (6.0,2.0);

					\foreach \i in {1,2,3,4,5,6} 
					{
						\draw(z) -- (x\i);
					}
					
					\draw(y1) -- (x1)--(y2);
					\draw(y3) -- (x2)--(y4);
					\draw(y5) -- (x3)--(y6);
					\draw(y7) -- (x4)--(y8);
					\draw(y9) -- (x5)--(y10);
					\draw(y11) -- (x6)--(y12);
					\draw  (y1) -- (y2) -- (y3) -- (y4) -- (y5) -- (y6);
					\draw(y7)-- (y8) -- (y9) -- (y10) -- (y11) -- (y12);

					\begin{scriptsize}
						\foreach \i in {1,2,3,4,5,6} 
						{
							\draw(x\i)[fill=white] circle(\vr);
						}
						\foreach \j in {1,2,3,4,5,6,7,8,9,10,11,12} 
						{
							\draw(y\j)[fill=white] circle(\vr);
						}
						\draw(z)[fill=white] circle(\vr);
						
					\end{scriptsize}
					\draw (0.5,5.5) node {$z$};
					\draw (0.5,2) node {$\ldots$};
					\draw (0.5,3) node {$\ldots$};
					\draw (-4.5, 3.5) node {$x_1$};
					\draw (5.5, 3.5) node {$x_k$};
					\draw (-5.3, 2.4) node {$y_1$};
					\draw (6.4, 2.4) node {$y_{2k}$};
				\end{scope}
				
			\end{tikzpicture}
			\caption{Graph $F_k$.
			}
			\label{fig:F_k}
		\end{center}
	\end{figure}

	For an integer $k$ with $k \ge 5$, define the graph $F_k$ as follows (see also Fig.~\ref{fig:F_k}). The vertex set is  $V(F_k)= \{z\} \cup V_x \cup V_y$
	with $V_x= \{x_i : i\in [k] \}$ and $V_y= \{y_i : i\in [2k] \}$, while the edge set is defined as
	$$E(F_k)= \{zx_i : i\in [k]\} \cup \{x_iy_{2i-1} : i\in [k] \} \cup \{x_iy_{2i} : i\in [k] \} \cup \{y_iy_{i+1} : i\in [2k-1]\}.$$
	\vspace{-4ex}
	
	\vspace{-3ex}
	
	\begin{proposition} \label{prop:F-k}
		For every integer $k\ge 5$,
		\begin{itemize}
			\item[(i)] $\mut(F_k)=2$, $\mud(F_k) =3$, $\muo(F_k)= \left\lceil \frac{2k}{3} \right\rceil$, and $\mu(F_k)\ge k+\left\lceil \frac{k}{2}\right\rceil$;
			\item[(ii)] $\mut(F_k-z)= \mud(F_k-z)= \muo(F_k-z)= \mu(F_k-z)= k+2$.
		\end{itemize}
	\end{proposition}
	\proof 
	\begin{description}
		\item[$\mut(F_k)=2$.]
		Clearly, $\{y_1, y_{2k}\}$ is a total mutual-visibility set of $F_k$. Thus $\mut(F_k) \ge 2$. Let $X$ be a $\mut$-set of $F_k$. 
		Observe that  $z \notin X$ since $z$ is the middle vertex of the convex $3$-path $x_1zx_2$. Since   $zx_iy_{2i}$ and $y_{j-1}y_{j} y_{j+1}$ are also convex $3$-paths in $F_k$, we infer that $x_i \notin X$, and $y_j \notin X$ hold for every $i \in [k]$  and $2 \le j \le 2k-1$. Hence, $X \subseteq \{y_1, y_{2k}\}$ and then $\mut(F_k) = 2$ follows.
		
		\item[$\mud(F_k) =3$.] 
		Since $F_k - \{x_1,y_1,y_2\}$ is an isometric subgraph of $F_k$, Lemma~\ref{lem:lower} implies that $\{x_1, y_1, y_2\}$ is a dual mutual-visibility set of $F_k$ and $\mud(F_k) \ge 3$. Let $X$ be a $\mud$-set of $F_k$. Since $k \ge 5$, there are two vertices $x_i$ and $x_j$ such that both are in $X$ or neither of them is in $X$. 
		However, if $z \in X$, then $x_i$ and $x_j$ are not $X$-visible, a contradiction. Thus $z \notin X$. 
		Consider now a vertex $x_i$ from $V_x$. It is the middle vertex of the convex $3$-paths 
		$y_{2i-1} x_iz$ and $y_{2i} x_iz$. Therefore, $x_i \in X$ is possible only if both $y_{2i-1} \in X$ and $y_{2i} \in X$ hold.
		Observe that $P^i: y_i \dots y_{i+3}$ is a convex $4$-path in $F_k$, for each $i \in [2k-3]$. By Lemma~\ref{lem:convex}, $X\cap V(P^i)$ is a dual mutual-visibility set in $P^i$. Therefore, if $y_{i+1} \in X$, then $X\cap V(P^i)= \{y_i, y_{i+1}\}$. Similarly, $y_{i+2} \in X$ implies $X\cap V(P^i)= \{y_{i+2}, y_{i+3}\}$. Since it holds for every $i \in [2k-3]$, we infer $X\cap V_y \subseteq \{y_1, y_2, y_{2k-1}, y_{2k}\}$. 
		
		Suppose now that $y_1, y_2 \in X$. Since $x_1y_2y_3$ is a convex $3$-path in $F_k$, vertices $x_1$ and $y_3$ are not $X$-visible. Since $X$ is a dual mutual-visibility set and $ y_3 \notin X$, we may infer that $x_1 \in X$. Then $y_1$ and $y_{2k-1}$ are not $X$-visible and hence  $y_{2k-1}\notin X$. Similarly, we get $y_{2k} \notin X$. In an analogous way, $y_{2k-1}, y_{2k} \in X$ implies $y_1, y_2 \notin X$.  We conclude $|X \cap V_y| \leq 2$. By our earlier observation, $x_i \notin X$ follows if $y_{2i-1} \notin X$ or $y_{2i} \notin X$ holds.
		
		To finish the proof, we consider the following cases. If $X \cap V_y= \emptyset$, then $X \cap V_x= \emptyset$ and $X=\emptyset$. If $|X \cap V_y|= 1$, then $X \cap V_x= \emptyset$ and $|X|=1$. Only the cases $X \cap V_y=\{y_1, y_2\}$ and $X \cap V_y=\{y_{2k-1}, y_{2k}\}$ allow $|X \cap V_x|= 1$. We conclude $|X| \leq 3$.  Therefore,  $\mud(F_k) = 3$.
		
		\item[$\muo(F_k)= \left\lceil \frac{2k}{3} \right\rceil$.] 
		Let $t=\left\lceil \frac{2k}{3} \right\rceil$. It is straightforward to see that $\{y_1, y_4, \ldots, y_{3t-2}\}$ is an outer mutual-visibility set of $F_k$. Hence, $\muo(F_k) \ge t$. Assume now that $X$ is a $\muo$-set of $F_k$. We will prove that $X$ contains at most $t$ vertices. 
		
		\noindent{\it Claim D.}\enskip If $y_i, y_j \in X$, then $|i-j| \ge 3$.
		
		\noindent{\it Proof.}\enskip Suppose, to the contrary, that $|i-j| \le 2$. Then we may choose an index $j'$ with $|i-j'|=3$ such that $y_j$ lies on the unique shortest $y_i, y_{j'}$-path. But in this case, $y_i$ and $y_{j'}$ are not $X$-visible and hence, $X$ is not a $\muo$-set, a contradiction. (If $j \in \{1,2k\}$ or $(i,j) \in \{(3,2), (2k-2,2k-1)\}$, we switch the roles of $y_i$ and $y_j$.) \smallqed
		
		\noindent{\it Claim E.}\enskip $|X \cap V_x| \leq 2$.
		
		\noindent{\it Proof.}\enskip Suppose for a contradiction that $\{x_i, x_j, x_\ell\} \subseteq X $ and $i <j < \ell$. Since $i+2 \le \ell$, the unique shortest $x_i, y_{2\ell}$-path goes through $x_\ell$. Thus $x_i$ and $y_{2\ell}$ are not $X$-visible, a contradiction.  \smallqed
		
		\noindent{\it Claim F.}\enskip If $z \in X$, then $|X|=1$.
		
		\noindent{\it Proof.}\enskip Suppose that $z \in X$. As $x_i$ lies on the unique shortest path $zx_iy_{2i}$, no vertex from $V_x$ belongs to $X$. Furthermore, as $z$ lies on the unique shortest $y_1,y_j$-path for every 
		$j \ge 6$ and on the unique shortest $y_j,y_{2k}$-path if $j \le 2k-5$, the present supposition $z \in X$ implies $y_j \notin X$ for these cases. Since $k \ge 5$, at least one of $j \ge 6$ and $j \le 2k-5$ holds for every $j \in [2k]$. Therefore, we have $|X|=1$ when $z \in X$. \smallqed
		
		\noindent{\it Claim G.}\enskip If $z \notin X$, then $|X| \leq t$.
		
		\noindent{\it Proof.}\enskip
		Suppose now that $x_i \in X$, for some $i\in [k]$. As $x_i$ lies on the unique shortest $y_j,y_{2i}$-path if $j \leq 2i-5$ and on the unique shortest $y_j,y_{2i-1}$-path if $2i+4 \leq j$, the $\muo$-set $X$ contains no vertex from $\{y_1, \dots , y_{2i-5}, y_{2i+4}, \dots ,y_{2k}\}$. As $x_iy_{2i-1}y_{2i-2}$ and $x_iy_{2i}y_{2i+1}$ are convex $3$-paths, $y_{2i-1} \notin X$ and $y_{2i} \notin X$ follow if $2 \leq i \leq k-1$. Consequently, 
		$$X \cap V_y \subseteq \{y_{2i-4}, y_{2i-3}, y_{2i-2}, y_{2i+1}, y_{2i+2}, y_{2i+3}\}
		$$
		holds whenever $x_i \in X$ for an integer $2 \leq i \leq k-1$. By Claim D, every two vertices in $X \cap V_y$ have a distance of at least $3$. Checking also the cases with $i \in \{1,k\}$, we may summarize our observations for the case of  $z \notin X$ as follows. 
		If $|X \cap V_y|\leq 2$, then Claim E implies $|X | \leq 4 \leq t$. If $|X \cap V_y|\ge 3$, then $|X \cap V_x| =0$. By Claim D, 
		$|X \cap V_y|\le t$, and we conclude  $|X| \le t$.  \smallqed
		
		Claims F and G imply $\muo(F_k) \leq  t$, and we may conclude the equality.
		
		\item[$\mu(F_k)\ge k+\left\lceil \frac{k}{2}\right\rceil$.] Observe that the following set $X$ is a mutual-visibility set in $F_k$, for every $k \ge 4$:
		$$X=\left\{x_{2i}: i \in \left\lfloor\frac{k}{2}\right\rfloor \right\} \cup 
		\left\{y_{4j-3}, y_{4j-2}: j \in \left\lceil\frac{k}{2}\right\rceil \right\}.
		$$
		Hence, $|X|= k+\left\lceil \frac{k}{2}\right\rceil$ verifies the statement. 
		
		\item[$\mut(F_k-z)= \mud(F_k-z)= \muo(F_k-z)= \mu(F_k-z)= k+2$].\\
		First, observe that $\{y_1, y_{2k}, x_1, x_2,  \ldots, x_k\}$ is a total mutual-visibility set in $F_k-z$ and therefore, $\mut(F_k -z) \ge k+2$. 
		We also show that $\mu(F_k-z) \le k+2$.
		Let $X$ be a $\mu$-set in $F_k-z$. 
		Since the path $P_{2k}: y_1\dots y_{2k}$ is a convex subgraph in $F_k-z$, Lemma~\ref{lem:convex}
		implies $|X \cap V(P_{2k})| \le 2$, and we may infer  $\mu(F_k-z) \le k+2$. Then, by~\eqref{eq:1}, we have $k+2 \le \mut(F_k-z) \le \mu(F_k-z) \le k+2$ and conclude $\mut(F_k-z)= \muo(F_k-z)= \mud(F_k-z) =\mu(F_k-z)=k+2$ as stated.
		\qed
	\end{description}
	It can be proved that $ k+\left\lceil \frac{k}{2}\right\rceil$ is, in fact, the exact value of $\mu(F_k)$. However, since the proof is not short and we do not use this equality in the continuation, the proof is omitted in this paper. Proposition~\ref{prop:F-k} directly implies the following statements.
	\begin{corollary} \enskip
		\label{cor:upper-vertex}
		\begin{itemize} 
			\item[(i)] There are no constants $a>0$ and $b$ such that $a\cdot \mud(G) +b \geq \mud(G-x)$ holds for every $G$ and $x \in V(G)$ when $G-x$ is connected. 
			\item[(ii)] There are no constants $a>0$ and $b$ such that $a\cdot \mut(G) +b \geq \mut(G-x)$ holds for every $G$ and $x \in V(G)$ when $G-x$ is connected. 
		\end{itemize}
	\end{corollary} 
	The values of the four visibility invariants typically can change drastically when a vertex is removed from the graph. However, we can point out that the increase in the mutual-visibility number is limited.
	\begin{theorem}
		\label{thm:mu-vertex}
		If $G$ is a graph, $x \in V(G)$, and $G-x$ is connected, then $\mu(G-x) \leq 2\, \mu(G)$ holds. 
	\end{theorem}
	\proof Let $G'$ denote $G-x$, and let $Y$ be a $\mu$-set in $G'$. Remark that $x \notin Y$. We will identify a set $B_Y \subseteq Y$, which is analogous to the ``blocking set'' defined in Section~\ref{sec:pre}. However, since $Y$ is not always a mutual-visibility set in $G$, we cannot use directly the definition and properties from Section~\ref{sec:pre}.
	
	For every $u \in Y$, fix a shortest $x,u$-path $P_{x,u}$ in $G$ such that $|V(P_{x,u} \cap Y|$ is minimum. 
	Then $P_{x,u}$ cannot have more than one internal vertex from $Y$. Suppose, to the contrary, that $P_{x,u}=x\dots u_1 \dots u_2 \dots u$ and $u_1, u_2 \in Y$. The subpath between $u_1$ and $u$ does not contain $x$ and therefore, it is a shortest $u_1,u$-path in $G'$. But then, since $Y$ is $\mu$-set in $G'$, there exists a shortest $u_1,u$-path that does not contain an internal vertex from $Y$. This contradicts the choice of $P_{x,u}$. 
	We say that a pair $(u',u)$ with $u',u \in Y$ is an \emph{$x$-blocking pair} in $G$ or that \emph{$u'$ $x$-blocks $u$}, if $P_{x,u}$ contains $u'$ and $u\neq u'$. Let $B_Y \subseteq Y$ contain a vertex  $u'\in Y$ if $u'$ $x$-blocks at least one vertex. We can prove, as before, that $|B_Y| \leq |Y|/2 $. 
	
	Finally, we show that $Y'= Y \setminus B_Y$ is a mutual-visibility set in $G$. Let $u$ and $v$ be two arbitrary vertices from $Y'$. If $d_G(u,v)=d_{G'}(u,v)$, then $u$ and $v$ are $Y'$-visible in $G$ via the path that ensures their $Y$-visibility in $G'$. Otherwise, $d_G(u,v) < d_{G'}(u,v)$ and every shortest $u,v$-path contains $x$ and therefore, $d_G(u,v)=d_G(u,x) + d_G(x,v)$. By concatenating $P_{x,u}$ and $P_{x,v}$, a shortest $u,v$-path $P$ is obtained. Since all $x$-blocking vertices were removed and $x \notin Y'$, no internal vertex of $P$ belongs to $Y'$. We conclude that $Y'$ is a mutual-visibility set in $G$ and $\mu(G) \ge |Y'| \ge |Y|/2 = \mu(G-x)/2$.
	\qed
	
	We do not have a sharp example for the inequality in Theorem~\ref{thm:mu-vertex}. First, we point out that such an example $G$ with a diameter of $2$ is impossible.  
	\begin{proposition}
		Let $G$ be a graph with $\diam(G)=2$, and let $x \in V(G)$ such that $G-x$ is connected. Then, $\mu(G-x) \leq \mu(G)$ holds.  
	\end{proposition}
	\proof Consider a $\mu$-set $Y$ in $G-x$ and two vertices $u,v \in Y$. If $u$ and $v$ are adjacent, they remain $Y$-visible in $G$. If $d_{G-x}(u,v)=2$, they remain $Y$-visible in $G$ via the same path as in $G-x$. If $d_{G-x}(u,v)>2$, then $d_{G}(u,v)=2$ and they are $Y$-visible via the path $uxv$. Therefore, $\mu(G) \ge |Y|=\mu(G-x)$.
	\qed

	Remark that $\mu(G-x) > \mu(G)$ is possible for graphs with $\diam(G)>2$. The graph $J$ constructed in Section~\ref{subsec:3.2} from a $12$-cycle and three leaves provides such an example.
	We have seen in Section~\ref{subsec:3.2} that $\mu(J)=3$. Deleting vertex $v_3$ from $J$, we get a tree with five leaves. Consequently, $\mu(J-v_3)=5= \frac{5}{3}\, \mu(J)$. Similar constructions can be obtained from longer cycles.
	
	\section{Open problems} \label{sec:open} 
	
	Only one remains open from the $16$ basic questions we considered (see Table~\ref{table:result}). We pose it as an open problem  for further consideration.
	\begin{problem} \label{prob:muo}
		Are there constants $a $ and $b$ such that $a\cdot \muo(G) +b \geq \muo(G-x)$ holds for every graph $G$ and $x \in V(G)$ when $G-x$ is connected?
	\end{problem}
	By Proposition~\ref{prop:F-k}, $\muo(F_k)= \left\lceil \frac{2k}{3}\right\rceil$ and $\muo(F_k-z) =k+2$ hold for every $k \ge 5$. Therefore, if a general upper bound $a\cdot \muo(G) +b$ on $\muo(G-x)$ exists, $a \ge \frac{3}{2}$ must be true.
	Indeed, if $k \ge \max\{5, b-2\}$, we obtain
	$$ a \ge \frac{\muo(F_k-z)-b}{\muo(F_k)} =\frac{k+2-b}{\left\lceil \frac{2k}{3}\right\rceil}\ge \frac{k+2-b}{\frac{2k+2}{3}}=\frac{3}{2} +\frac{3-3b}{2k+2}\, 
	$$
	and, for every $\epsilon >0$, integer $k$ can be chosen such that $\left| \frac{3-3b}{2k+2} \right| < \epsilon$.
	\medskip

	Concerning the sharpness of the lower and upper bounds stated in Theorems~\ref{thm:edge-mu}, \ref{thm:edge-outer}, and~\ref{thm:mu-vertex}, several questions remain open. 
	\begin{problem} \label{prob:2}
		Is it possible to improve the upper bound in Theorem~\ref{thm:edge-mu}, the inequalities in Theorem~\ref{thm:edge-outer}, and the upper bound in Theorem~\ref{thm:mu-vertex}?
	\end{problem}
	We proved that the lower bound in Theorem~\ref{thm:edge-mu} is asymptotically sharp that is, a general lower bound $a\, \mu(G) +b \leq \mu(G-e)$ with a constant $a >1/2$ is impossible. However, a slight improvement might be achieved. We conjecture that $\frac{1}{2} \mu(G) +1 \leq \mu(G-e)$ might hold for every graph. If this is true,  the sharpness is demonstrated by the graphs $H_k$ (see Proposition~\ref{prop:Hk}). 
	\section*{Acknowledgments}
	
	This work was supported by the Slovenian Research and Innovation Agency (ARIS) under the grants P1-0297 and N1-0355.
	
	\baselineskip11pt

\end{document}